# $n$−absorbing $I$−primary ideals in commutative rings

**Sarbast A. Anjuman[1], Ismael Akray[2]**
[1,2]Mathematics deparment, Soran University, Soran, Kurdistan region-Erbil, Iraq
*Corresponding author: Sarbast A.Anjuman
sarbast.mohammed@soran.edu.iq

**Abstract**

We define a new generalization of $n$−absorbing ideals in commutative rings called $n$−absorbing $I$−primary ideals. We investigate some characterizations and properties of such new generalization. If $P$ is an $n$−absorbing $I$−primary ideal of $R$ and $\sqrt{IP} = I\sqrt{P}$, then $\sqrt{P}$ is a $n$−absorbing $I$−primary ideal of $R$. And if $\sqrt{P}$ is an $(n-1)$−absorbing ideal of $R$ such that $\sqrt{I\sqrt{P}} \subseteq IP$, then $P$ is an $n$−absorbing $I$−primary ideal of $R$.

## 1 Introduction

In our article all rings are commutative ring with non-zero identity. In the recent years many generalizations of prime ideals were defined. Here state some of them. The notion of a weakly prime ideal was introduced by Anderson and Smith, where a proper ideal $P$ of a commutative ring $R$ is a weakly prime if $x, y \in R$, and $0 \neq xy \in P$ then $x \in P$ or $y \in P$ in [1] . Ebrahimi Atani and Farzalipour defined the nation of weakly primary ideals [2] . The authors in [3] and [4] introduced the notions $2$−absorbing and $n$−absorbing ideals in commutative rings. A proper ideal $P$ is said to be $2$−absorbing (or $n$−absorbing) ideal if whenever the product of three (or $n+1$) elements of $R$ in $P$, the product of two (or $n$) of these elements is in $P$.

In [5] and [6], the author Akray introduced the notions $I$−prime ideal, $I$−primary ideal and $n$−absorbing $I$−ideal in commutative rings as a generalization of prime ideals. For fixed proper ideal $I$ of a commutative ring $R$ with identity, a proper ideal $P$ of $R$ is an $I$−primary if for $c, d \in R$ with $cd \in P - IP$, then $c \in P$ or $d \in \sqrt{P}$. Throughout the paper the notation $b_1 \cdots \widehat{b_i} \cdots b_n$ means that $b_i$ is excluded from the product $b_1 \cdots b_n$. A proper ideal $P$ of $R$ is an $n$−absorbing $I$−primary ideal if for $b_1, \cdots, b_{n+1} \in R$ such that $b_1 \cdots b_{n+1} \in P - IP$, then $b_1 \cdots b_n \in P$ or $b_1 \cdots b_{i-1} b_{i+1} \cdots b_{n+1} \in \sqrt{P}$ for some $i \in \{1,2,\cdots,n\}$ where $\sqrt{P}$ is the radical of the ideal $P$.

Assume that $R$ is an integral domain with quotient field $F$. The authors in [7] introduced a "proper ideal $P$ of $R$ is a strongly primary if, whenever $cd \in P$ with $c, d \in F$, we have $c \in P$ or $d \in \sqrt{P}$. In [8], a proper ideal $P$ of $R$ is a strongly $I$−primary ideal if $cd \in P - IP$ with $c, d \in F$, then $c \in P$ or $d \in \sqrt{P}$. It is said that a proper ideal $P$ of $R$ is quotient $n$−absorbing $I$−primary" if $b_1 b_2 \cdots b_{n+1} \in P$ with $b_1, b_2, \ldots, b_{n+1} \in F$, then $b_1 b_2 \cdots b_n \in P$ or $b_1 \cdots \widehat{b_i} \cdots b_{n+1} \in \sqrt{P}$ for some $1 \leq i \leq n$. Set $P$ is an ideal of a ring $R$, let $P$ be an $n$−absorbing $I$−primary ideal of $R$ and $b_1, \ldots, b_{n+1} \in R$. The statement is that $(b_1, \ldots, b_{n+1})$ is an $I - (n+1)$−tuple of $P$ if $b_1 \cdots b_{n+1} \in IP$, $b_1 b_2 \cdots b_n \notin P$ and for any $1 \leq i \leq n$, $b_1 \cdots \widehat{b_i} \cdots b_{n+1} \notin \sqrt{P}$.

# 2  $n-$absorbing $I-$primary ideals

In this section, we start with to define the definition of an $n-$absorbing $I-$primary ideal of a ring $R$.

**Definition**  *A proper ideal $P$ of $R$ is an $I-$primary if for $c,d \in R$ with $cd \in P - IP$, then $c \in P$ or $d \in \sqrt{P}$.*

**Definition**  *A proper ideal $P$ of $R$ is an $n-$absorbing $I-$primary ideal if for $b_1, \cdots, b_{n+1} \in R$ such that $b_1 \cdots b_{n+1} \in P - IP$, then $b_1 \cdots b_n \in P$ or $b_1 \cdots b_{i-1} b_{i+1} \cdots b_{n+1} \in \sqrt{P}$ for some $i \in \{1, 2, \cdots, n\}$.*

**Example 2.1**  *Consider the ring $A = k[t_1, t_2, \ldots, t_{n+2}]$, where $k$ is a field and suppose that $P = \langle t_1 t_2 \cdots t_{n+1}, t_1^2 t_2 \cdots t_n, t_1^2 t_{n+2} \rangle$, $I = \langle t_1 t_2 \cdots t_n, t_1 t_2 \cdots t_{n+1} \rangle$. Then $P - IP = \langle t_1 \cdots t_{n+1}, t_1^2 t_2 \cdots t_n, t_1^2 t_{n+2} \rangle - \langle t_1 t_2 \cdots t_{n+1}, t_1^2 \cdots t_n, t_1^2 \cdots t_{n+1}, t_1^2 t_2 \cdots t_{n+2} \rangle$. Hence $P$ is an $n-$absorbing $I-$primary ideal but $P$ is not $n-$absorbing.*

**Proposition 2.2**  *We set that $R$ is a ring. Based on this, the following statements can be considered equivalent:*

*(i) $P$ is an $n-$absorbing $I-$primary ideal of $R$;*

*(ii) For any elements $\alpha_1, \ldots, \alpha_n \in R$ with $\alpha_1 \cdots \alpha_n$ not in $\sqrt{P}$, $(P:_R \alpha_1 \cdots \alpha_n) \subseteq \left[ \cup_{i=1}^{n-1} \left( \sqrt{P}:_R \alpha_1 \cdots \hat{\alpha}_i \cdots \alpha_n \right) \right] \cup (P:_R \alpha_1 \cdots \alpha_{n-1}) \cup (IP:_R \alpha_1 \cdots \alpha_n)$.*

*Proof.* (i) $\Rightarrow$ (ii) Set $\alpha_1, \ldots, \alpha_n \in R$ such that $\alpha_1 \cdots \alpha_n \notin \sqrt{P}$. Let $r \in (P:_R \alpha_1 \cdots \alpha_n)$. So $r\alpha_1 \cdots \alpha_n \in P$. If $r\alpha_1 \cdots \alpha_n \in IP$, then $r \in (IP:_R \alpha_1 \cdots \alpha_n)$. Let $r\alpha_1 \cdots \alpha_n \notin IP$. Since $\alpha_1 \cdots \alpha_n \notin \sqrt{P}$, either $r\alpha_1 \cdots \alpha_{n-1} \in P$, that is, $r \in (P:_R \alpha_1 \cdots \alpha_{n-1})$ or for some $1 \leq i \leq n-1$ we have $r\alpha_1 \cdots \hat{\alpha}_i \cdots \alpha_n \in \sqrt{P}$, that is, $r \in \left( \sqrt{P}:_R \alpha_1 \cdots \hat{\alpha}_i \cdots \alpha_n \right)$ or $\alpha_1 \alpha_2 \cdots \alpha_n \in \sqrt{P}$. By assumption the last case is not hold. Consequently
$(P:_R \alpha_1 \cdots \alpha_n) \subseteq \left[ \cup_{i=1}^{n-1} \left( \sqrt{P}:_R \alpha_1 \cdots \hat{\alpha}_i \cdots \alpha_n \right) \right] \cup (P:_R \alpha_1 \cdots \alpha_{n-1}) \cup (IP:_R \alpha_1 \cdots \alpha_n)$.

(ii) $\Rightarrow$ (i) Let $\beta_1 \beta_2 \cdots \beta_{n+1} \in P - IP$ for some $\beta_1, \beta_2, \ldots, \beta_{n+1} \in R$ such that $\beta_1 \beta_2 \cdots \beta_n \notin P$. Then $\beta_1 \in (P:_R \beta_2 \cdots \beta_{n+1})$. If $\beta_2 \cdots \beta_{n+1} \in \sqrt{P}$, then we are done. Hence we may set $\beta_2 \cdots \beta_{n+1} \notin \sqrt{P}$ and so by part (ii), $(P:_R \beta_2 \cdots \beta_{n+1}) \subseteq \left[ \cup_{i=2}^{n} \left( \sqrt{P}:_R \beta_2 \cdots \hat{\beta}_i \cdots \beta_{n+1} \right) \right] \cup (P:_R \beta_2 \cdots \beta_n) \cup (IP:_R \beta_2 \cdots \beta_{n+1})$. Since $\beta_1 \beta_2 \cdots \beta_{n+1} \notin IP$ and $\beta_1 \beta_2 \cdots \beta_n \notin P$, the only possibility is that $\beta_1 \in \cup_{i=2}^{n} \left( \sqrt{P}:_R \beta_2 \cdots \hat{\beta}_i \cdots \beta_{n+1} \right)$. Then $\beta_1 \beta_2 \cdots \hat{\beta}_i \cdots \beta_{n+1} \in \sqrt{P}$ for some $2 \leq i \leq n$. Hence $P$ is an $n-$absorbing $I-$primary ideal of $R$. □

**Proposition 2.3** *If $V$ be a valuation domain with the quotient field $F$. Then all $n$-absorbing $I$-primary ideal of $V$ is a quotient $n$-absorbing $I$-primary ideal of $R$.*

*Proof.* We can certainly assume that $P$ is $n$-absorbing $I$-primary ideal of $V$, and $a_1 a_2 \cdots a_{n+1} \in P$ for some $a_1, a_2, \ldots, a_{n+1} \in V$ such that $a_1 a_2 \cdots a_n \notin P$. If $a_{n+1} \notin V$, then $a_{n+1}^{-1} \in V$. So $a_1 \cdots a_n a_{n+1} a_{n+1}^{-1} = a_1 \cdots a_n \in P$, which is a contradiction. So $a_{n+1} \in V$. If $a_i \in V$ for all $1 \leq i \leq n$, then there is nothing to prove. If $a_i \notin V$ for some $1 \leq i \leq n$, then $a_1 \cdots \widehat{a_i} \cdots a_{n+1} \in P \subseteq \sqrt{P}$. Consequently, $P$ is a quotient $n$-absorbing $I$-primary. □

**Proposition 2.4** *Set $P$ be an $n$-absorbing $I$-primary ideal of $R$ such that $\sqrt{IP} = I\sqrt{P}$, then $\sqrt{P}$ is a $n$-absorbing $I$-primary ideal of $R$.*

*Proof.* Let us assume $a_1 a_2 \cdots a_{n+1} \in \sqrt{P} - I\sqrt{P}$ for some $a_1, a_2, \ldots, a_{n+1} \in R$ such that $a_1 \cdots \widehat{a_i} \cdots a_{n+1} \notin \sqrt{P}$ for every $1 \leq i \leq n$. Thus we have $n \in \mathbb{N}$ such that $a_1^n a_2^n \cdots a_{n+1}^n \in P$. If $a_1^n a_2^n \cdots a_{n+1}^n \in IP$, then $a_1 a_2 \cdots a_{n+1} \in \sqrt{IP} = I\sqrt{P}$, which is a contradiction. Since $P$ is an $n$-absorbing $I$-primary, our hypothesis implies $a_1^n a_2^n \cdots a_n^n \in P$. So $a_1 a_2 \cdots a_n \in \sqrt{P}$ and $\sqrt{P}$ is an $n$-absorbing $I$-primary ideal of $R$. □

**Theorem 2.5** *Assume that for any $1 \leq i \leq k$, $I_i$ is an $n_i$-absorbing $I$-primary ideal of $R$ such that $\sqrt{P_i} = q_i$ is an $n_i$-absorbing $I$-primary ideal of $R$, respectively. Let $n = n_1 + n_2 + \cdots + n_k$. The following statements does hold:*

(1) $P_1 \cap P_2 \cap \cdots \cap P_k$ *is an $n$-absorbing $I$-primary ideal of $R$.*
(2) $P_1 P_2 \cdots P_k$ *is an $n$-absorbing $I$-primary ideal of $R$.*

*Proof.* The proof of the two parts are similar, so we prove just the first. Let $H = P_1 \cap P_2 \cap \cdots \cap P_k$. Then $\sqrt{H} = P_1 \cap P_2 \cap \cdots \cap P_k$. Let $a_1 a_2 \cdots a_{n+1} \in H - IH$ for some $a_1, a_2, \ldots, a_{n+1} \in R$ and $a_1 \cdots \widehat{a_i} \cdots a_{n+1} \notin \sqrt{H}$ for any $1 \leq i \leq n$. By, $\sqrt{H} = P_1 \cap P_2 \cap \cdots \cap P_k$ is an $n$-absorbing $I$-primary, then $a_1 a_2 \cdots a_n \in P_1 \cap P_2 \cap \cdots \cap P_k$. We prove that $a_1 a_2 \cdots a_n \in H$. For all $1 \leq i \leq k$, $P_i$ is an $n_i$-absorbing $I$-primary and $a_1 a_2 \cdots a_n \in P_i - IP_i$, then we have $1 \leq \beta_1^i, \beta_2^i, \ldots, \beta_{n_i}^i \leq n$ such that $a_{\beta_1^i} a_{\beta_2^i} \cdots a_{\beta_{n_i}^i} \in P_i$. If $\beta_r^l = \beta_s^m$ it is for two couples $l, r$ and $m, s$, then

$$a_{\beta_1^1} a_{\beta_2^1} \cdots a_{\beta_{n_1}^1} \cdots a_{\beta_1^l} a_{\beta_2^l} \cdots a_{\beta_r^l} \cdots a_{\beta_{n_l}^t} \cdots$$

$$a_{\beta_1^m} a_{\beta_2^m} \cdots \widehat{a_{\beta_m^m}} \cdots a_{\beta_{n_m}^m} \cdots a_{\beta_1^k} a_{\beta_2^k} \cdots a_{\beta_{n_k}^k} \in \sqrt{H}$$

Therefore $a_1 \cdots \widehat{a_{\beta_m^m}} \cdots a_n a_{n+1} \in \sqrt{H}$, which is a contradiction. So $\beta_j^i$ 's are distinct. Hence $\{a_{\beta_1^1}, a_{\beta_2^1}, \ldots, a_{\beta_{n_1}^1}, a_{\beta_1^2}, a_{\beta_2^2}, \ldots, a_{\beta_{n_2}^2}, \ldots, a_{\beta_1^k}, a_{\beta_2^k}, \ldots, a_{\beta_{n_k}^k}\} = \{a_1, a_2, \ldots, a_n\}$. If $a_{\beta_1^i} a_{\beta_2^i} \cdots a_{\beta_{n_i}^i} \in P_i$ for any $1 \leq i \leq k$, then

$$a_1 a_2 \cdots a_n = a_{\beta_1^1} a_{\beta_2^1} \cdots a_{\beta_{n_1}^1} a_{\beta_1^2} a_{\beta_2^2} \cdots a_{\beta_{n_2}^2} \cdots a_{\beta_1^k} a_{\beta_2^k} \cdots a_{\beta_{n_k}^k} \in H,$$

thus we are done. Therefore we may assume that $a_{\beta_1^1} a_{\beta_2^1} \cdots a_{\beta_{n_1}^1} \notin P_1$. Since $P_1$ is $I =$

$n_1$ −absorbing $I$ −primary and
$$a_{\beta_1^1}a_{\beta_2^1}\cdots a_{\beta_{n_1}}a_{\beta_1^2}a_{\beta_2^2}\cdots a_{\beta_{n_2}^2}\cdots a_{\beta_1^k}a_{\beta_2^k}\cdots a_{\beta_{n_k}^k}a_{n+1}=a_1\cdots a_{n+1}\in P_1-IP_1,$$
then we have $a_{\beta_1^2}a_{\beta_2^2}\cdots a_{\beta_{n_2}^2}\cdots a_{\beta_1^k}a_{\beta_2^k}\cdots a_{\beta_{n_k}^k}a_{n+1}\in P_1$. On the other hand $a_{\beta_1^2}a_{\beta_2^2}\cdots a_{\beta_{n_2}^2}\cdots a_{\beta_1^k}a_{\beta_2^k}\cdots a_{\beta_{n_k}^k}a_{n+1}\in P_2\cap\cdots\cap P_k$. Consequently $a_{\beta_1^2}a_{\beta_2^2}\cdots a_{\beta_{n_2}^2}$ $\cdots a_{\beta_1^k}a_{\beta_2^k}\cdots a_{\beta_{n_k}^k}a_{n+1}\in\sqrt{H}$, which is a contradiction. Similarly $a_{\beta_1^i}a_{\beta_2^i}\cdots a_{\beta_{n_i}^i}\in P_i$ for every $2\leq i\leq k$. Then $a_1a_2\cdots a_n\in H$. □

**Proposition 2.6** *Assume that $P$ is an ideal of a ring $R$ with $\sqrt{I\sqrt{P}}\subseteq IP$. If $\sqrt{P}$ is an $(n-1)$ −absorbing ideal of $R$, then $P$ is an $n$ −absorbing $I$ −primary ideal of $R$.*

*Proof.* Let $\sqrt{P}$ be an $(n-1)$ −absorbing, and consider $b_1b_2\cdots b_{n+1}\in P-IP$ for some $b_1,b_2,\ldots,b_{n+1}\in R$ and $b_1b_2\cdots b_n\notin P$. Since
$$(b_1b_{n+1})(b_2b_{n+1})\cdots(b_nb_{n+1})=(b_1b_2\cdots b_n)b_{n+1}^n\in P\subseteq\sqrt{P}-I\sqrt{P}.$$
Then for some $1\leq i\leq n$,
$$(b_1b_{n+1})\cdots\widehat{(x_ib_{n+1})}\cdots(b_nb_{n+1})=(b_1\cdots\hat{b}_i\cdots b_n)b_{n+1}^{n-1}\in\sqrt{P},$$
and so $b_1\cdots\hat{b}_i\cdots b_nb_{n+1}\in\sqrt{P}$. Consequently $P$ is an $n$ −absorbing $I$ −primary ideal of $R$. □

We recall that a proper ideal $Q$ of $R$ is an **$n$ −absorbing primary** if $a_1,a_2,\cdots,a_{n+1}\in R$ and $a_1a_2\cdots a_{n+1}\in Q$, then $a_1a_2\cdots a_n\in Q$ or the product of $a_{n+1}$ with $(n-1)$ of $a_1,a_2,\cdots,a_n$ is in $\sqrt{Q}$. It is clearly every $n$ −absorbing primary is an $n$ −absorbing I −primary.

**Proposition 2.7** *Suppose that $R$ is a ring and $r\in R$, a nonunit and $m\geq 2$ is not negative integer. Let $(0:_R r)\subseteq\langle a\rangle$, then $\langle r\rangle$ is an $n$ −absorbing $I$ −primary, for some $I$ with $IP\subseteq I^m$ if and only if $\langle a\rangle$ is $n$ −absorbing primary.*

*Proof.* Let $\langle r\rangle$ be $n$ −absorbing $I^m$ −primary, and $a_1a_2\cdots a_{n+1}\in\langle r\rangle$ for some $a_1,a_2,\ldots,a_{n+1}\in R$. If $a_1a_2\cdots a_{n+1}\notin\langle r^m\rangle$, then $a_1a_2\cdots a_n\in\langle r\rangle$ or $a_1\cdots\hat{a}_i\cdots a_{n+1}\in\sqrt{\langle r\rangle}$ for some $1\leq i\leq n$. Based on this assumption, $a_1a_2\cdots a_{n+1}\in\langle r^m\rangle$. Hence $a_1a_2\cdots a_n(a_{n+1}+r)\in\langle r\rangle$. If $a_1a_2\cdots a_n(a_{n+1}+r)\notin\langle r^m\rangle$, then $a_1a_2\cdots a_n\in\langle r\rangle$ or $a_1\cdots\hat{a}_i\cdots a_n(a_{n+1}+r)\in\sqrt{\langle r\rangle}$ for some $1\leq i\leq n$. So $a_1a_2\cdots a_n\in\langle r\rangle$ or $a_1\cdots\hat{a}_i\cdots a_{n+1}\in\sqrt{\langle r\rangle}$ for some $1\leq i\leq n$. Hence, suppose that $a_1a_2\cdots a_n(a_{n+1}+r)\in\langle r^m\rangle$. Thus $a_1a_2\cdots a_{n+1}\in\langle r^m\rangle$ implies that $a_1a_2\cdots a_nr\in\langle r^m\rangle$. Therefore, there exists $s\in R$ such that $a_1a_2\cdots a_n-sr^{m-1}\in(0:_R r)\subseteq\langle r\rangle$. Consequently $a_1a_2\cdots a_n\in\langle r\rangle$.

□

**Proposition 2.8** *Assume $V$ is a valuation domain and $n \in \mathbb{N}$. Let $P$ be an ideal of $V$ such that $P^{n+1}$ is not principal. Then $P$ is an $n-$absorbing $I^{n+1}-$primary if and only if it is an $n-$absorbing primary.*

*Proof.* ($\Rightarrow$) Let $P$ be an $n-$absorbing $I^n-$primary that is n't $n-$absorbing primary. Therefore there are $a_1, \ldots, a_{n+1} \in R$ such that $a_1 \cdots a_{n+1} \in P$, but neither $a_1 \cdots a_n \in P$ nor $a_1 \cdots \widehat{a_i} \cdots a_{n+1} \in \sqrt{P}$ for any $1 \leq i \leq n$. Hence $\langle a_i \rangle \nsubseteq P$ for any $1 \leq i \leq n+1$. And so $V$ is a valuation domain, thus $P \subset \langle a_i \rangle$ for any $1 \leq i \leq n+1$, and so $P^{n+1} \subseteq \langle a_1 \cdots a_{n+1} \rangle$. Therefore $P^{n+1}$ is not principal, then $a_1 \cdots a_{n+1} \in P - P^{n+1}$. Therefore $P$ is an $n-$absorbing $I^{n+1}-$primary implies that either $a_1 \cdots a_n \in P$ or $a_1 \cdots \widehat{a_i} \cdots a_{n+1} \in \sqrt{P}$ for some $1 \leq i \leq n$, which is a contradiction. Hence $P$ is an $n-$absorbing primary ideal of $R$.
($\Leftarrow$) Is trivial.  □

**Theorem 2.9** *We consider that $J \subseteq P$ are a proper ideals of a ring $R$.*

(1) *Let $P$ is an $n-$absorbing $I-$primary ideal of $R$, then $P/J$ is a $n-$absorbing $I-$primary ideal of $R/J$.*

(2) *Let $J \subseteq IP$ and $P/J$ be an $n-$absorbing $I-$primary ideal of $R/J$, then $P$ is an $n-$absorbing $I-$primary ideal of $R$.*

(3) *Let $IP \subseteq J$ and $P$ be an $n-$absorbing $I-$primary ideal of $R$, then $P/J$ is a weakly $n-$absorbing primary ideal of $R/J$.*

(4) *Let $JP \subseteq IP$, $J$ be an $n-$absorbing $I-$primary ideal of $R$ and $P/J$ be a weakly $n-$absorbing primary ideal of $R/J$, then $P$ is an $n-$absorbing $I-$primary ideal of $R$.*

*Proof.* (1) Set $b_1, b_2, \ldots, b_{n+1} \in R$ such that $(b_1 + J)(b_2 + J) \cdots (b_{n+1} + J) \in (P/J) - I(P/J) = (P/J) - (I(P) + J)/J$. Then $b_1 b_2 \cdots b_{n+1} \in P - IP$ and from being $P$ is an $n-$absorbing $I-$primary, we obtain $b_1 \cdots b_n \in P$ or $b_1 \cdots \widehat{b_i} \cdots b_{n+1} \in \sqrt{P}$ for some " $1 \leq i \leq n$. And so $(b_1 + J) \cdots (b_n + J) \in P/J$ or $(b_1 + J) \cdots (\widehat{b_i + J}) \cdots (b_{n+1} + J) \in \sqrt{P}/J = \sqrt{P/J}$ for some $1 \leq i \leq n$. Hence we prove that $P/J$ is $n-$absorbing $I-$primary ideal of $R/J$."

(2) Set $b_1 b_2 \cdots b_{n+1} \in P - IP$ for some $b_1, b_2, \ldots, b_{n+1} \in R$. Then $(b_1 + J)(b_2 + J) \cdots (b_{n+1} + J) \in (P/J) - (I(P)/J) = (P/J) - I(P/J)$. From being $P/J$ is an $n-$absorbing $I-$primary, we obtain that $(b_1 + J) \cdots (b_n + J) \in P/J$ or $(b_1 + J) \cdots (\widehat{b_i + J}) \cdots (b_{n+1} + J) \in \sqrt{P/J} = \sqrt{P}/J$ for some $1 \leq i \leq n$. Therefore $b_1 \cdots b_n \in P$ or $b_1 \cdots \widehat{b_i} \cdots b_{n+1} \in \sqrt{P}$ for some $1 \leq i \leq n$, hence $P$ is an $n-$absorbing $I-$primary ideal of $R$.

(3) Resulted directly from part (1).

(4) Set $b_1 \cdots b_{n+1} \in P - IP$ where $b_1, \ldots, b_{n+1} \in R$. Note that $b_1 \cdots b_{n+1} \notin JP$ because $JP \subseteq IP$. If $b_1 \cdots b_{n+1} \in J$, then either $b_1 \cdots b_n \in J \subseteq P$ or $b_1 \cdots \widehat{b_i} \cdots b_{n+1} \in \sqrt{J} \subseteq \sqrt{P}$ for some $1 \leq i \leq n$, since $J$ is an $n-$absorbing $I-$primary. If $b_1 \cdots b_{n+1} \notin J$, then $(b_1 + J) \cdots (b_{n+1} + J) \in (P/J) - \{0\}$ and so either $(b_1 + J) \cdots (b_n + J) \in P/J$ or $(b_1 + J) \cdots (\widehat{b_i + J}) \cdots (b_{n+1} + J) \in \sqrt{P/J} = \sqrt{P}/J$ for some $1 \leq i \leq n$. Therefore $b_1 \cdots b_n \in P$ or $b_1 \cdots \widehat{b_i} \cdots b_{n+1} \in \sqrt{P}$ for some $1 \leq i \leq n$. Hence $P$ is an $n-$absorbing $I-$primary ideal of $R$.

**Proposition 2.10** *Suppose that $P$ is an ideal of a ring $R$ such that $IP$ is an $n-$absorbing primary ideal of $R$. If $P$ is an $n-$absorbing $I-$primary ideal of $R$, then $P$ is an $n-$absorbing primary ideal of $R$.*

*Proof.* Let $a_1 a_2 \cdots a_{n+1} \in P$ for some elements $a_1, a_2, \ldots, a_{n+1} \in R$ such that $a_1 a_2 \cdots a_n \notin P$. If $a_1 a_2 \cdots a_{n+1} \in IP$, then $IP$ $n-$absorbing primary and $a_1 a_2 \cdots a_n \notin IP$ implies that $a_1 \cdots \widehat{a_i} \cdots a_{n+1} \in \sqrt{IP} \subseteq \sqrt{P}$ for some $1 \leq i \leq n$, and so we are done. When $a_1 a_2 \cdots a_{n+1} \notin IP$ clearly the result follows.

**Theorem 2.11** *If $P$ is an $n-$absorbing $I-$primary ideal of a ring $R$ and $(a_1, \ldots, a_{n+1})$ is an $I-(n+1)-$tuple of $P$ for some $a_1, \ldots, a_{n+1} \in R$. Then for every elements $\alpha_1, \alpha_2, \ldots, \alpha_m \in \{1, 2, \ldots, n+1\}$ which $1 \leq m \leq n$,*
$$a_1 \cdots \widehat{a_{\alpha_1}} \cdots \widehat{a_{\alpha_2}} \cdots \widehat{a_{\alpha_m}} \cdots a_{n+1} I^m \subseteq IP$$

*Proof.* We claim that by using induction on $m$. We take $m = 1$ and assume $a_1 \cdots \widehat{a_{\alpha_1}} \cdots a_{n+1} x \notin IP$ for some $x \in P$. Then $a_1 \cdots \widehat{a_{\alpha_1}} \cdots a_{n+1}(a_{\alpha_1} + x) \notin IP$. Since $P$ is a $n-$absorbing $I-$primary ideal of $R$ and $a_1 \cdots \widehat{a_{\alpha_1}} \cdots a_{n+1} \notin P$, we conclude that $a_1 \cdots \widehat{a_{\alpha_1}} \cdots \widehat{a_{\alpha_2}} \cdots a_{n+1}(a_{\alpha_1} + x) \in \sqrt{P}$, for some $1 \leq \alpha_2 \leq n+1$ different from $\alpha_1$. Hence $a_1 \cdots \widehat{a_{\alpha_2}} \cdots a_{n+1} \in \sqrt{P}$, a contradiction. Thus $a_1 \cdots \widehat{a_{\alpha_1}} \cdots a_{n+1} P \subseteq IP$. Here assume that $m > 1$ and for every integers less than $m$ the prove does hold. Let $a_1 \cdots \widehat{a_{\alpha_1}} \cdots \widehat{a_{\alpha_2}} \cdots \widehat{a_{\alpha_m}} \cdots a_{n+1} x_1 x_2 \cdots x_m \notin IP$ for some $x_1, x_2, \ldots, x_m \in P$. According to the induction assumption, we conclude that there exists $\zeta \in IP$ such that
$$a_1 \cdots \widehat{a_{\alpha_1}} \cdots \widehat{a_{\alpha_2}} \cdots \widehat{a_{\alpha_m}} \cdots a_{n+1}(a_{\alpha_1} + x_1)(a_{\alpha_2} + x_2) \cdots (a_{\alpha_m} + x_m)$$
$$= \zeta + a_1 \cdots \widehat{a_{\alpha_1}} \cdots \widehat{a_{\alpha_2}} \cdots \widehat{a_{\alpha_m}} \cdots a_{n+1} x_1 x_2 \cdots x_m \notin IP.$$
Now, we have two cases.

Case 1. Set $\alpha_m < n+1$. Since from being $P$ is an $n-$absorbing $I-$primary, then
$$a_1 \cdots \widehat{a_{\alpha_1}} \cdots \widehat{a_{\alpha_2}} \cdots \widehat{a_{\alpha_m}} \cdots a_n(a_{\alpha_1} + x_1)(a_{\alpha_2} + x_2) \cdots (a_{\alpha_m} + x_m) \in P,$$
or
$$a_1 \cdots \widehat{a_{\alpha_1}} \cdots \widehat{a_{\alpha_2}} \cdots \widehat{a_{\alpha_m}} \cdots \widehat{a_j} \cdots a_{n+1}(a_{\alpha_1} + x_1)(a_{\alpha_2} + x_2) \cdots (a_{\alpha_m} + x_m)$$
$$\in \sqrt{P}$$
for some $j < n+1$ distinct from $\alpha_i$'s; or
$$a_1 \cdots \widehat{a_{\alpha_1}} \cdots \widehat{a_{\alpha_2}} \cdots \widehat{a_{\alpha_m}} \cdots a_{n+1}(a_{\alpha_1} + x_1) \cdots \widehat{(a_{\alpha_i} + x_i)} \cdots (a_{\alpha_m} + x_m) \in \sqrt{P}$$
for some $1 \leq i \leq m$. Thus either $a_1 a_2 \cdots a_n \in P$ or $a_1 \cdots \widehat{a_j} \cdots a_{n+1} \in \sqrt{P}$ or $a_1 \cdots \widehat{a_{\alpha_i}} \cdots a_{n+1} \in \sqrt{P}$, which each of these cases which is a contradiction.

Case 2. Set $\alpha_m = n+1$. Since from being $P$ is an $n-$absorbing $I-$primary, then
$a_1 \cdots \widehat{a_{\alpha_1}} \cdots \widehat{a_{\alpha_2}} \cdots \widehat{a_{\alpha_{m-1}}} \cdots \widehat{a_{n+1}}(a_{\alpha_1} + x_1)(a_{\alpha_2} + x_2) \cdots (a_{\alpha_m} + x_m) \in P$, or
$$a_1 \cdots \widehat{a_{\alpha_1}} \cdots \widehat{a_{\alpha_2}} \cdots \widehat{a_{\alpha_{m-1}}} \cdots \widehat{a_j} \cdots \widehat{a_{n+1}}(a_{\alpha_1} + x_1)(a_{\alpha_2} + x_2) \cdots (a_{\alpha_m} + x_m)$$
$$\in \sqrt{P},$$
for some $j < n+1$ different from $\alpha_i$'s; or $a_1 \cdots \widehat{a_{\alpha_1}} \cdots \widehat{a_{\alpha_2}} \cdots \widehat{a_{\alpha_{m-1}}} \cdots \widehat{a_{n+1}}(a_{\alpha_1} + x_1) \cdots (\widehat{a_{\alpha_i} + x_i}) \cdots (a_{\alpha_m} + x_m) \in \sqrt{P}$ for some $1 \leq i \leq m-1$. Thus either $a_1 a_2 \cdots a_n \in P$ or $a_1 \cdots \widehat{a_j} \cdots a_{n+1} \in \sqrt{P}$ or $a_1 \cdots \widehat{a_{\alpha_i}} \cdots a_{n+1} \in \sqrt{P}$, which each of these cases which are a contradiction. Thus
$$a_1 \cdots \widehat{a_{\alpha_1}} \cdots \widehat{a_{\alpha_2}} \cdots \widehat{a_{\alpha_m}} \cdots a_{n+1} I^m \subseteq IP$$

**Theorem 2.12** *If $P$ is an $n-$absorbing $I-$primary ideal of $R$ which is not an $n$-absorbing primary ideal. Then*
  (i) $P^{n+1} \subseteq IP$.
  (ii) $\sqrt{P} = \sqrt{IP}$.

*Proof.* (i) Since $P$ is assumed not to be an $n-$absorbing primary ideal of $R$, so $P$ has an $I-(n+1)-$tuple zero $(b_1, \ldots, b_{n+1})$ for some $b_1, \ldots, b_{n+1} \in R$. Let $c_1 c_2 \cdots c_{n+1} \notin IP$ for some $c_1, c_2, \ldots, c_{n+1} \in P$. Therefore, according to the Theorem 2.11, there is $\lambda \in IP$ such that $(b_1 + c_1) \cdots (b_{n+1} + c_{n+1}) = \lambda + c_1 c_2 \cdots c_{n+1} \notin IP$. Hence either $(b_1 + c_1) \cdots (b_n + c_n) \in P$ or $(b_1 + c_1) \cdots (\widehat{b_i + c_i}) \cdots (b_{n+1} + c_{n+1}) \in \sqrt{P}$ for some $1 \leq i \leq n$. Thus either $b_1 \cdots b_n \in P$ or $b_1 \cdots \widehat{b_i} \cdots b_{n+1} \in \sqrt{P}$ for some $1 \leq i \leq n$, which is a contradiction. Hence $P^{n+1} \subseteq IP$.
  (ii) Clearly, $\sqrt{IP} \subseteq \sqrt{P}$. As $P^{n+1} \subseteq IP$, we obtain $\sqrt{P} \subseteq \sqrt{IP}$, we are done.